\documentclass{amsart}
\usepackage{amssymb,  latexsym, ulem}
\if0
\usepackage{mathtools}
\mathtoolsset{showonlyrefs=true/fault} %このパッケージは`mathtools.sty`とケンカする場合がある．それを避けるために
\usepackage{showlabels} %`[inline]` option (see https://konoyonohana.blog.fc2.com/blog-entry-299.html)
\fi

\usepackage{cancel}

\usepackage[all]{xy}

\newtheorem{theorem}{Theorem}[section]
\newtheorem{lemma}[theorem]{Lemma}
\newtheorem{sublemma}[theorem]{Sublemma}
\newtheorem{corollary}[theorem]{Corollary}
\newtheorem{proposition}[theorem]{Proposition}

\theoremstyle{definition}
\newtheorem{definition}[theorem]{Definition}

\theoremstyle{remark}
\newtheorem{remark}[theorem]{Remark}

\newcommand{\lem}[2]{\begin{lemma}\label{#1} #2\end{lemma}}

\newcommand{\prop}[2]{\begin{proposition}\label{#1}#2\end{proposition}}

\newcommand{\thm}[2]{\begin{theorem}\label{#1}#2\end{theorem}}

\newcommand{\cor}[2]{\begin{corollary}\label{#1}#2\end{corollary}}

\newcounter{th}

0

\newcommand{\al}{\alpha}
\newcommand{\be}{\beta}
\newcommand{\ga}{\gamma}

\newcommand{\de}{\delta}

\newcommand{\ep}{\varepsilon}

\newcommand{\et}{\eta}

\newcommand{\io}{\iota}

\newcommand{\la}{\lambda}
\newcommand{\La}{{\Lambda}}

\newcommand{\Si}{{\Sigma}}

\newcommand{\ph}{\varphi}

\newcommand{\ps}{\psi}

\newcommand{\Om}{{\Omega}}
\newcommand{\barr}{\begin{array}}
\newcommand{\ear}{\end{array}}
\newcommand{\bsk}{\begin{array}{rcl}}
\newcommand{\esk}{\end{array}}
\newcommand{\skh}{\begin{eqnarray*}}
\newcommand{\sko}{\end{eqnarray*}}
\newcommand{\bcss}{\begin{cases}}
\newcommand{\ecss}{\end{cases}}

\newcommand{\AL}[1]{\begin{align*}#1\end{align*}}

\newcommand{\AR}[2]{$$\begin{array}{#1}#2\end{array}$$}
\newcommand{\ARr}[3]{$$\begin{array}{#1}#2\end{array}\lnr{#3}$$}
\newcommand{\aR}[2]{$\begin{array}{#1}#2\end{array}$}

\newcommand{\cass}[1]{\begin{cases}#1\ecss}
\newcommand{\ak}{\quad}

\newcommand{\XY}[1]{$$\xymatrix@=15pt{#1}$$}
\newcommand{\XYr}[2]{$$\vcenter{\xymatrix@=15pt{#1}}\lnr{#2}$$}
\renewcommand{\Xy}[2]{$$\xymatrix@=#1pt{#2}$$}
\newcommand{\Xyr}[3]{$$\vcenter{\xymatrix@=#1pt{#2}}\lnr{#3}$$}
\newcommand{\Xyrs}[4]{$$\vcenter{\xymatrix@=#1pt{#2}}\lnrs{#4}{#3}$$}
\newcommand{\xR}[2]{$\xymatrix@R#1pt{#2}$}
\newcommand{\xyR}[2]{$$\xymatrix@R#1pt{#2}$$}
\newcommand{\xyT}[1]{$$\xymatrix@R10pt{#1}$$}
\newcommand{\xyRr}[3]{$$\vcenter{\xymatrix@R#1pt{#2}}\lnr{#3}$$}
\newcommand{\xyRrs}[4]{$$\vcenter{\xymatrix@R#1pt{#2}}\lnrs{#4}{#3}$$}

\newcommand{\arr}{\ar[r]}

\newcommand{\ard}{\ar[d]}

\newcommand{\wt}[1]{\widetilde{#1}}

\newcommand{\Z}{\mathbb Z}

\def\o+{\oplus}
\newcommand{\Op}{\bigoplus}
\newcommand{\x}{\times}
\newcommand{\ox}{\otimes}

\newcommand{\e}{{\rm Ext}}
\newcommand{\lnb}{\refstepcounter{theorem}\leqno(\thetheorem)}
\newcommand{\lnr}[1]{\lnb\label{#1}}
\newcommand{\nr}{\refstepcounter{theorem}\thetheorem}
\newcommand{\kko}[1]{(\ref{#1})}
\newcommand{\mbx}[1]{\quad\mbox{#1}\quad}
\newcommand{\cf}{{\it cf.}\ }
\newcommand{\sus}{\subset }
\newcommand{\sm}{\wedge }
\newcommand{\Ker}{{\rm Ker}\ }

\newcommand{\im}{{\rm Im}\, }
\newcommand{\Hom}{\mbox{\rm Hom}}

\newcommand{\cg}{\equiv}
\newcommand{\lrk}[1]{\left\langle #1\right\rangle}

\newcommand{\qand}{\mbx{and}}

\newcommand{\Lt}{\left}
\newcommand{\Rt}{\right}
\newcommand{\LR}[1]{\Lt(#1\Rt)}
\newcommand{\cln}{\colon}
\renewcommand{\O}[1]{\overline{#1}}

\newcommand{\ANSS}{Adams-Novikov spectral sequence}
\newcommand{\mx}[1]{\begin{matrix}#1\end{matrix}}

\renewcommand{\dh}{\mathfrak{h}}
\newcommand{\dH}{\mathfrak{H}}
\newcommand{\di}{\mathfrak{i}}

\newcommand{\F}{{\mathbb F}_p}
\newcommand{\xar}{\xrightarrow}

\newcommand{\cL}{{\mathcal L}}

\newcommand{\cS}{{\mathcal S_{(p)}}}

\newcommand{\Vee}{\bigvee}

\newcommand{\p}{\noindent {\it Proof.} }
\newcommand{\q}{\hfill $\qed$ \medskip}

\newcommand{\h}[2]{h_{#1,#2}}

\renewcommand{\b}[2]{b_{#1,#2}}

\newcommand{\thick}[1]{\text{\rm thick}\lrk{#1}}

\newcommand{\rank}{\mbox{\rm rank }}

\newcommand{\der}{\partial}

\newcommand{\mdi}[2]{\medskip
\noindent
(\nr)  \label{#1} {\it #2 }
\medskip
}

\newcommand{\madi}[1]{\noindent (\nr)\label{#1}}

\newcommand{\Ex}{\mbox{\rm Ex}}

\begin{document}

%\title{A note on $\be_1$-action on the stable homotopy of spheres }
\title[A note on the Greek letter elements]{A note on the Greek letter elements of the homotopy of the $L_n$-local spheres}
\author{Ryo Kato}
\address{Department of Core Studies, Kochi University of Technology, Kami,
Kochi, 782-8502, Japan}
\email{ryo\_kato\_1128@yahoo.co.jp}
\author{ Katsumi Shimomura}
\address{Department of Mathematics, Faculty of Science, Kochi University, Kochi,
780-8520, Japan}
\email{katsumi@kochi-u.ac.jp}
\author{Mao-no-suke Shimomura}
\address{Department of Mathematics, Faculty of Science and Technology, Kochi University, Kochi, 780-8520, Japan}
\email{mao.shimomura@outlook.jp}
%\date{\today}

%\subjclass[2010]{Primary~55Q45, Secondary~55T15, 55Q51}
%\keywords{stable homotopy groups of spheres, Greek letter elements}

\begin{abstract}
\if0
We consider the homotopy groups of the $v_3^{-1}BP$-local (resp. $v_4^{-1}BP$) sphere spectrum at prime five (resp. seven), and show the groups contains the gamma (resp. delta) family.  
\fi
Let $BP$ denote the Brown-Peterson spectrum at a prime $p$, whose homotopy groups are isomorphic to the polynomial algebra generated by elements $v_i$'s for $i\ge 1$.
We consider the homotopy groups of the $v_n^{-1}BP$-localized sphere spectrum $L_nS^0$ with $n^2\le 2p-1$, and show that the groups contain the $n$-th Greek letter family.
For the proof of this, we further show the existence of the $v_n^{-1}BP$-localized Smith-Toda spectrum $W_n$ for the case $n^2\le 2p-1$.
If the Smith-Toda spectrum $V(n-1)$ exists, then $L_nV(n-1)$ is an example of $W_n$.
Previously,
$W_n$ is shown to exist if $n^2+n\le 2p-1$. 

We also consider the case where $(p,n)=(7,4)$, and show the existence of the delta family in $\pi_*(L_4S^0)$.  
\end{abstract}
\maketitle

\section{Introduction}

Let $\cS$ denote the stable homotopy category of the $p$-local spectra for a prime number $p$.
In $\cS$, we have the Brown-Peterson spectrum $BP\in \cS$ with homotopy groups $\pi_*(BP)=BP_*=\Z_{(p)}[v_1,v_2,\dots\,]$ for $v_n\in \pi_{2(p^n-1)}(BP)=BP_{2(p^n-1)}$.
Then  we obtain  a Bousfield localization functor $L_n\cln \cS\to \cS$ with respect to $v_n^{-1}BP$ along with the natural transformation $\et_n\cln id \to L_n$.
%Here, $v_n\in \pi_{2(p^n-1)}(BP)=BP_{2(p^n-1)}$ denotes an indecomposable generator.

For studying the homotopy groups $\pi_*(X)$ of a spectrum $X$, 
we have the Adams-Novikov spectral sequence
$$
E_2^{s,t}(X)=\e_{BP_*BP}^{s,t}(BP_*,BP_*(X))\Longrightarrow \pi_{t-s}(X).
$$
In the $E_2$-term $E_2^{s,t}(S^0)$ of this spectral sequence for $\pi_*(S^0)$,
we have the generators called  the $n$-th Greek letter elements $\al_s^{(n)}\in E_2^{n,t}(S^0)$ for $t=2s(p^n-1)-2\sum_{k=1}^{n-1} (p^k-1)$.
It is well known that $\al^{(n)}_s$ is a permanent cycle of the spectral sequence in the following cases:
\ARr {ll}{
n=1:&\al_s(=\al_s^{(1)}) \mbx{is a permanent cycle if $p\ge 3$ and $s\ge 1$,}\\ 
&\hspace{1in}\mbox{ or if $p=2$ and $1\le s\not\cg 3$ $(4)$,}\\
n=2:&\be_s(=\al_s^{(2)}) \mbx{is a permanent cycle if $p\ge 5$ and $s\ge 1$,}\\ 
&\hspace{1in}\mbox{ or if $p=3$ and $1\le s\not\cg 4,7,8$ $(9)$,}\\
n=3:&\ga_s(=\al_s^{(3)}) \mbx{is a permanent cycle if $p\ge 7$ and $s\ge 1$.}
}{known}
Let $E(n)(\in\cS)$ denote the $n$-th Johnson-Wilson spectrum. Then we have the $E(n)$-based Adams spectral sequence
$$
E_2^{s,t}(X)_n=\e_{E(n)_*E(n)}^{s,t}(E(n)_*,E(n)_*(X))\Longrightarrow \pi_{t-s}(L_nX), \lnr{ANSS}
$$
and we have a canonical map $\ell_n\cln BP\to E(n)$, which induces a homomorphism $(\ell_n)_*\cln BP_*=\Z_{(p)}[v_1,v_2,\dots\,] \to E(n)_*=\Z_{(p)}[v_1,\dots, v_n,v_n^{-1}]$ sending $v_k$ to $v_k$ for $k\le n$ and to zero for $k>n$. 
The map $\ell_n$ induces a map 
$(\ell_n)_*\cln \{E_r^{s,t}(X)\} \to \{E_r^{s,t}(X)_n\}$ of spectral sequences from the Adams-Novikov spectral sequence to the $E(n)$-based Adams spectral sequence. 
By the map, the fact \kko{known} implies that
\AR {ll}{
p=2,\ n=1:& \Z/p\lrk{\al_s} \sus \pi_{2s-1}(L_1S^0) \mbx{unless $s\cg 3$ $(4)$,}\\
p=3,\ n=2:& \Z/p\lrk{\be_s}\sus \pi_{16s-6}(L_2S^0) \mbx{unless $s\cg 4,7,8$ $(9)$.}
} 
In this paper, we consider the Greek letter elements in $\pi_*(L_nS^0)$.
First we have  results  for   $(p,n)=(5,3)$ and $(7,4)$.

\thm{main}{For each integer $s\ne 0$,
 the third Greek letter element $\ga_s\in E_2^3(S^0)_3$ at $p=5$ and the fourth Greek letter element $\de_s\in E_2^4(S^0)_4$ at $p=7$
are permanent cycles in the respective spectral sequences. These permanent cycles detect generators of submodules of $\pi_*(L_3S^0)$ and $\pi_*(L_4S^0)$ isomorphic to $\Z/p$.
}

The case for $(p,n)=(5,3)$ is generalized as follows:

\thm{main2}{Let $p$ and $n$ satisfy $n^2\le 2p-1$ and $(p,n)\ne (2,1)$, $(3,2)$, and $s$ be a non-zero integer.
Then, the $n$-th Greek letter element $\al^{(n)}_s\in E_2^n(S^0)_n$ at $p$
is a permanent cycle, which yields a generator of a submodule of $\pi_*(L_nS^0)$ isomorphic to $\Z/p$.
}

Theorems \ref{main} and \ref{main2} imply the following:
\cor{main2c}{
We have the fourth Greek letter element $\de_s\in \pi_*(L_4S^0)$ of order $p$  if $p\ge 7$.
}

Theorem \ref{main2} is based on the existence of the $L_n$-local Smith-Toda spectrum $W_n$, which is defined by the $BP_*$-homology 
$$
BP_*(W_n)=v_n^{-1}BP_*/I_n
$$
 for the invariant prime ideal $I_n=(p,v_1,\dots , v_{n-1})$ of $E(n)_*$.
 We notice that $W_n$ is not determined uniquely.
Since $E(n)_*(X)=E(n)_*\ox_{BP_*}BP_*(X)$, we also have
$$
E(n)_*(W_n)=E(n)_*/I_n
$$
We notice that  if the Smith-Toda spectrum $V(n-1)$ exists, the localized spectrum $L_nV(n-1)$ is one of $W_n$'s.
The theorem \cite[Th.~5.7]{sy} implies that if  $n^2+n\le 2p-1$, then $W_n$ exists uniquely.
Consider a condition on $(p,n)$:
$$
p\cg 3\ (4),\ \text{ or \ \ $n$ is even.}\lnr{pn}
$$
\if0
Let $MJ$ denote the generalized Moore spectrum such that $BP_*(MJ)=BP_*/J$ for an invariant regular ideal $J$.
Then, $V(n-1)=MI_n$.
On the existence, we read off from \cite[Th.~5.7]{sy} that $L_nMJ_n$ for an invariant regular ideal $J_n$ of length $n$ exists if $n^2+n\le 2p-1$.
\fi
Then, we improve the results:

\thm{main3}{ Let $n^2\le  2p-1$ and $(p,n)\ne (2,1)$, $(3,2)$.
Then, there exists an $L_n$-local Smith-Toda spectrum $W_n$ at the prime $p$.
Furthermore, if \kko{pn} holds, then it is unique.
%
%Furthermore, if $n^2\le 2p-2$, then there exists $L_nMJ_n$ for an invariant regular ideal $J_n\sus E(n)_*$ of length $n$.
}

In the next section, we compute some cohomology $H^*S(n)$ of the Morava stabilizer algebra $S(n)$.
Section 3 explains the existence of the $L_n$-local Smith-Toda spectrum $W_n$ based on \cite{sy}, and we show Theorem \ref{main3}.
We study the dual map $D_n(i_n)\cln D_n(W_n)\to S^0$ to define Greek letter elements in section 4.
 In the last section, we show Theorems \ref{main} and \ref{main2}.

\if0

In the homotopy groups $\pi_*(S^0)$ of the sphere spectrum $S^0$, we have the $n$-th Greek letter elements $\al_s^{(n)}$.
Some of these elements survives to $\pi_*(L_nS^0)$ under the homomorphism $(\et_n)_*\cln \pi_*(S^0)\to \pi_*(L_nS^0)$.

and consider a full subcategory $\cL_n$ of $\cS$ consisting of $v_n^{-1}BP$-local spectra.

Let $p$ be an odd prime number.

The Greek letter elements are the generators of the stable homotopy groups of spheres at a prime $p$.
It is known that $\al_s$, $\be_s$ and $\ga_s$ exists if $p\ge 3$, $\ge 5$ and $\ge 7$, respectively.
At $p=2$ and $p=3$, $\al_s$ and $\be_s$ exist for a restricted $s$:
$\al_{4t+3}\not\in\pi_*(S^0)$ at the prime two, and $\be_{9t+4}$, $\be_{9t+7}$, $\be_{9t+8}\not\in\pi_*(S^0)$.
 
 $d_3(v_1^3/2)=v_1h_0^3/2$ in $E_3^*(L_1S)$ and 
 
For $(p,n)$ with $n<p-1$ and $n^2\le 2p-2$, we have $n$-th Greek letter 
family $\{\al^{(n)}_s\mid s\ge 1\}$ 
in $\pi_*(L_nS)$. 

In this paper, we obtain the following theorem:
\thm{1}{
Let $(p,n)=(5,3)$, $(7,4)$.
Then, we have a gamma family $\{\ga_s\mid s\in\Z\}\sus \pi_*(L_3S)$ at the prime five
and a delta family $\{\de_s\mid s\in\Z\}\sus \pi_*(L_4S)$ at the prime seven.
\if0
 be a positive integer with $n<p-1$.
If a ring spectrum $L_nV(n-1)$ exists and $n^2\le 4p-4$, then there exists a family $\{\al^{(n)}_s\mid s\ge 1\}$ of Greek letter elements in $\pi_*(L_nS)$. 
\fi
}

\fi

\section{The cohomology groups $H^{*,*}S(n)$}

In this section, we show $H^{s,t}S(n)=0$ for some $s$ and $t$ including $H^{q+1,q}S(n)=0$ for $n^2\le q+1=2p-1$.
We also prove  $H^{q+1,q}S(n)=0$ for  $(p,n)=(7,4)$.
Hereafter,  we set
$$
q=2p-2.
$$
We notice that the Morava vanishing theorem says that $H^{s,*}S(n)=0$ if $(p-1)\nmid n$ and $s> n^2$.
Furthermore, we see that
\begin{itemize}
\item $n^2\le 2p-1\implies n<p$, \  and so $$(p-1)\mid n\iff (p,n)=(2,1), (3,2).$$
 %n=p-1. %\lnr{n2}
\item $n^2\le 2p-1\ \text{and \kko{pn} holds}\implies  n^2\le 2p-2$.
\end{itemize}
For a positive integer $n$, the relation $n^2\le 2p-1$ implies $n<p$, and so $(p-1)\mid n$ if and only if $n =p-1$.
If $n=p-1$, then $p^2-2p+1\le 2p-1$, $0\ge p^2-4p+2=(p-2)^2-2$ and $p=2$ or 3.

We consider a condition on $(p,n)$:
$$
n^2\le 2p-1\qand n\le p-2. \lnr{cond}
$$
Under $n^2\le 2p-1$, 
$$
n\le p-2\iff (p,n)\ne (2,1), (3,2)\iff (p-1)\nmid n.
$$

Our computation is based on results of Ravenel's:\\

\madi{ravenel} (\cite[6.3.4.~Th., 6.3.5~Th.]{r:book})   If $0\le n\le p-2$, then there are spectral sequences
\begin{enumerate}
\item
$E_2=H^*(L(n,n))\Rightarrow H^*(E_0S(n)^*)$,
\item
$E_2=H^*(E_0S(n)^*)\Rightarrow H^*S(n)$.

\end{enumerate}

\medskip

\madi{r638} (\cite[6.3.8~Th.]{r:book})
$H^*(L(n,n))$ is the cohomology of the exterior complex $E(\h ij\mid 1\le i\le n$ and $j\in \Z/n)$, with differential 
$$
d(\h ij)=\sum_{s=1}^{i-1}\h sj\h{i-s}{s+j}.
$$
Here, the degree $|\h ij|\in \Z/|v_n|$ of $\h ij$ is given by
$$
|\h ij|=2p^j(p^i-1)=p^je(i)q\mbx{for}  e(i)=\sum_{k=0}^{i-1}p^k=\frac{p^i-1}{p-1}.
$$

Under the condition \kko{cond}, the spectral sequences in \kko{ravenel} collapse from the $E_2$-term, and
these imply the following lemma:

\lem{HS}{Suppose \kko{cond}. Then,
$H^{*}E(\h ij)_n\cong H^*S(n)$ as a module.
Here, $E(\h ij)_n$ denotes the complex of \kko{r638}.
In particular, $H^sS(n)=0$ if $s>n^2$.
}

We consider a lexicographic ordered set $\dH_n=\{(i,j)\mid 1\le i\le n, j\in \Z/n\}$.
Put $h_A=\prod_{(i,j)\in A}\h ij$ and $h_A^*=h_{\dH_n\setminus A}=\prod_{(i,j)\in \dH_n\setminus A}\h ij$ for a subset $A\sus \dH_n$. 
Then,  $E(\h ij)_n$ is additively generated by $h_A$ for $A\sus \dH_n$, and there is an isomorphism
$$
-^*\colon E(\h ij)^{s,t}_n\cong E(\h ij)_n^{n^2-s, -t} \lnr{isoD}
$$
assigning $h_A^*$ to $h_A$. Indeed, $h_A^*h_A=\pm g_n$ for $g_n=h_{\dH_n}$.

\if0 H^{s,s}S(n)=0:

\lem{00}{Suppose \kko{cond}. %that $n^2\le 2p-1$ and $(p,n)\ne (2,1), (3,2)$. 
Then, %$H^{n^2,*}(S(n))=0$ and
$H^{s,s}(S(n))=0$ for $s>0$.
}

\p
It suffices to show $E(\h ij)_n^{s,s}=0$ for $s>0$ by Lemma \ref{HS}.
By degree reason, $E(\h ij)_n^{s,s}=0$ unless $q\mid s$.
For $s=kq$, $E(\h ij)_n^{kq,kq}=0$ if $k\ge 2$ or $k=1$ and $n^2< 2p-2=q$.
If $n^2=q$ (resp. $=q+1$), then  by \kko{isoD}, $E(\h ij)_n^{q,q}\cong E(\h ij)_n^{0,-q}=0$ (resp. $E(\h ij)_n^{q,q}\cong E(\h ij)_n^{1,-q}=0$).
\q
\fi

\lem{zero}{
%$H^{q+1,q}S(n)=0$ for $n^2\le 2p-1$ with $(p,n)\ne (2,1), (3,2)$, and 
$H^{13,12}S(4)=0$ at $p=7$. 
}
\p
\if0
If $n^2<q+1$ and $(p-1)\nmid n$, then $H^{q+1,q}S(n)=0$ by the Morava vanishing theorem (\cf \cite[1.4.11]{r:book}).
Suppose that $n^2=q+1$.
Then, $H^{q+1,*}L(n,n)\cong \Z/p\{g_n\}$ by \kko{r638}, where $g_n=\prod_{0\le i<n, j\in\Z/n}\h ij$.
It follows that $H^{q+1,q}L(n,n)= 0$ by degree reason, which implies $H^{q+1,q}S(n)=0$ by \kko{ravenel}.
\fi
%Now turn to the case for $(p,n)=(7,4)$. Then, 
In this case, 
$q=12$ and $|v_4|=4800=12\x 400$.
Since $E(\h ij)^{13,12}_4\cong E(\h ij)^{3,-12}_4$ by \kko{isoD}, we first consider $E(\h ij)^{3,-12}_4$.
An additive generators of $E(\h ij)^{3,-12}_4$ has a form of one of
$$
\h ij \h4a\h4b, \ak \h ij \h kl \h4a\qand \h ij \h k\ell \h mn,
$$
in which $i$, $k$ and $m$ are less than four.
We have a table of elements and its degrees:
\AR{|c||c|c|c|c||c|c|c|c|}{
\hline
\text{elements}&\h10&\h11&\h12&\h13&\h20&\h21&\h22&\h23\\
\hline
\text{degree}/12&1&7&49&-57&8&56&-8&-56\\
\hline
\cline{1-5}
\text{elements}&\h30&\h31&\h32&\h33\\
\cline{1-5}
\text{degree}/12&57&-1&-7&-49\\
\cline{1-5}
}
and $|\h 4a|=0$.
\begin{description}
\item[
The elements of the form $\h ij\h4a\h4b$]
The degree of the elements are read off from the above table, and we obtain
$$
\h31\h4a\h4b.
$$
for $a,b\in\Z/4$.
\item[
The elements of the form $\h ij\h k\ell\h4a$]
From the table, we obtain the following elements of degree $-12$: %=399\x 12$:
$$
\h11\h22\h4a\qand \h13\h21\h4a
$$
for $a\in\Z/4$.
\item[
The elements of the form $\h ij\h k\ell\h mn$] We assume that $i\le k\le m$.
If $|\h i0\h k\ell\h mn|/12=a$, then $|\h ij\h k{\ell+j}\h m{n+j}|/12=7^ja$.
We consider elements $\h ij\h k{\ell+j}\h m{n+j}\in E(\h ij)_4^{3,-12}$. In this case, the degree gives rise to
$7^ja\cg -1$ mod $(400)$. In $\Z/(400)$, we deduce $7^{-1}=-57$, $7^{-2}%=57\x (7^2+7+1) =(-1)\x (7+1)+57
=49$, $7^{-3}%=-57\x 49
=7$ from $7^4=1$.
Thus, an element $\h i0\h k\ell\h mn$ with $a=-1, 57, -49,-7$ yields an element $\h ij\h k{\ell+j}\h m{n+j}$ of desired degree.
So we will find elements $\h k\ell\h mn$ of degree $12b$ for
$b=a-c$ for $c=1,8$ and $57$ if $i=1,2$ and $3$, respectively, and
$$
b=\cass{-2,56,-50,-8&i=1\\-9,49,-57,-15 &i=2\\ -58,0,-106,-64&i=3}
$$
%
%$=-2,56,-50,-8$ .
From the above table, we read off the following elements:

\smallskip

\noindent
\aR{c}{
\|\h11\h12\|=56, \ak \|\h11\h13\|=-50,\ak \|\h12\h13\|=-8,\ak \|\h30\h31\|=56,\\
 \|\h31\h32\|=-8,\ak \|\h31\h33\|=-50\ak\mbx{(if $i=1$),}\\
\|\h22\h31\|=-9,\ak \|\h21\h32\|=49,\ak \|\h22\h30\|=49,\ak \|\h22\h33\|=-57,\\
\|\h23\h31\|=-57,\ak \|\h22\h32\|=-15\ak\mbx{(if $i=2$).}
}

\if0

If $|\h10\h k\ell\h mn|/12=a$, then \linebreak $|\h1j\h k{\ell+j}\h m{n+j}|/12=7^ja$.
We consider elements $\h1j\h k{\ell+j}\h m{n+j}\in E(\h ij)^{3,-12}$. In this case, the degree gives rise to
$7^ja\cg -1$ mod $(400)$. In $\Z/(400)$, we deduce $7^{-1}=-57$, $7^{-2}%=57\x (7^2+7+1) =(-1)\x (7+1)+57
=49$, $7^{-3}%=-57\x 49
=7$ from $7^4=1$.
Thus, an element $\h10\h k\ell\h mn$ with $a=-1, 57, -49,-7$ yields an element of desired degree.
So we will find elements $\h k\ell\h mn$ of degree $12b$ for
$b=a-1=-2,56,-50,-8$.
From the above table, we read off the following elements:

\smallskip

\aR{c}{
\|\h11\h12\|=56, \ak \|\h11\h13\|=-50,\ak \|\h12\h13\|=-8,\ak \|\h30\h31\|=56,\\
 \|\h31\h32\|=-8,\ak \|\h31\h33\|=-50,
}
\fi
\smallskip

\noindent
in which $\|x\|=|x|/12$.
We notice that  no element arises for $i=3$.
These give us the elements

\smallskip

\noindent
$
\h11\h12\h13, \ \ \h12\h13\h11,\ \  \h13\h11\h12,\ \  \h11\h31\h32,\ \ \h13\h30\h31,\ \  \h12\h33\h31,
$  

\noindent
$
\h20\h22\h31, \ \ \h21\h22\h33,\ \  \h21\h23\h31,\ \  \h22\h20\h31,\ \ \h22\h21\h33,\ \  \h23\h21\h31.
$  

\smallskip

\end{description}

The above calculation shows that $E(\h ij)^{13,12}_4$ is generated by
\ARr l{
a_{k,\ell}=(\h31\h4k\h4\ell)^*, \ak b_k=(\h13\h21\h4k)^*, \ak b'_k=(\h11\h22\h4k)^*, \\
c_0= (\h11\h12\h13)^*, \ak
c_1= (\h11\h31\h32)^*, \ak c_2=(\h12\h31\h33)^*, \\ c_3=(\h13\h30\h31)^*,\ak 
c_4=(\h20\h22\h31)^*, \ak c_5=(\h21\h22\h33)^*, \\ c_6=(\h21\h23\h31)^*
 }{genE}
 (see \kko{isoD} for $x^*$).
% Here, $x^*$ denotes an element such that $xx^*=\pm g_n$ for the generator $g_n\in H^{q+1,*}L(n,n)$ given at the beginning of the proof.
 The differential $d$ on the exterior complex $E(\h ij)_4$ acts on these elements as follows in which $x_k=(\h31\h4k)^*$:
 \AL{
 d((\h11\h22\h4k\h4\ell)^*)&=(\h31\h4k\h4\ell)^*=a_{k,\ell}\\
 d((\h11\h12\h13\h4k)^*)&=(\h13\h21\h4k)^*+(\h11\h22\h4k)^*=b_k+b_k'\\
 d((\h13\h21\h4k)^*)&=d(b_k)=(\h31\h4k)^*=x_k\\
 d((\h11\h12\h13)^*)&=d(c_0)=(\h13\h21)^*+(\h11\h22)^*\\
%}
 %\AL{
 d((\h11\h31\h32)^*)&=d(c_1)=(\h31\h41)^*+(\h31\h42)^*=x_1+x_2\\
 d((\h12\h31\h33)^*)&=d(c_2)=(\h31\h42)^*+(\h31\h43)^*=x_2+x_3\\
 d((\h13\h31\h30)^*)&=d(c_3)=(\h31\h43)^*+(\h31\h40)^*=x_3+x_0\\
 d((\h20\h22\h31)^*)&=d(c_4)=(\h31\h40)^*+(\h31\h42)^*=x_0+x_2\\
d((\h21\h22\h33)^*)&=d(c_5)=0 %\\
}
\AL{d((\h21\h23\h31)^*)&=d(c_6)=(\h31\h41)^*+(\h31\h43)^*=x_1+x_3\\
%}
%\AL{
d((\h11\h13\h21\h32)^*)&=(\h13\h21\h41)^*+(\h13\h21\h42)^*+ (\h11\h31\h32)^*\\
&=b_1+b_2+c_1\\
d((\h12\h13\h20\h31)^*)&=(\h22\h20\h31)^*+(\h13\h30\h31)^*+(\h12\h33\h31)^*\\
&=c_4+c_3+c_2\\
d((\h11\h13\h22\h30)^*)&=(\h13\h31\h30)^*+(\h11\h22\h43)^*+(\h11\h22\h40)^*\\
&=c_3+b_3'+b_0'\\
%d((\h10\h13\h21\h31)^*)&=(\h23\h21\h31)^*+(\h13\h30\h31)^*+(\h13\h21\h40)^*+(\h13\h21\h41)^*\\
d((\h10\h11\h22\h31)^*)&=(\h20\h22\h31)^*+(\h11\h32\h31)^*+(\h11\h22\h41)^*+(\h11\h22\h40)^*\\
&=c_4+c_1+b_1'+b_0'\\
d((\h11\h12\h22\h33)^*)&=(\h21\h22\h33)^*+(\h12\h31\h33)^*+(\h11\h22\h42)^*+(\h11\h22\h43)^*\\
&=c_5+c_2+b_2'+b_3'\\
d((\h11\h12\h23\h31)^*)&=(\h21\h23\h31)^*+(\h12\h33\h31)^*+(\h11\h32\h31)^*\\
 &=c_6+c_2+c_1
 }
whose coefficients are up to sign.
Then, $\im d$ contains $A=\{a_{k,\ell}\mid 0\le k<\ell\le 3\}$ and $B=\{b_k\pm b'_k\mid 0\le k\le 3\}$, and $c_0$ is not a cocycle.
For each $c_i$ $(1\le i\le 6)$, we have a cocycle $\wt c_i=c_i+\la_jb_j+\la_kb_k$ for some $0\le j,k\le 3$ and $\la_\ell\in\{-1,0,1\}$ by the third equality.
We see that $C=\{\wt c_i\mid 1\le i\le 6\}\sus\im d$ by the last six equalities together with the second equality.
\if0
The first four relations show that the generators $(\h31\h4k\h4\ell)^*$, $(\h13\h21\h4k)^*$, $(\h11\h22\h4k)^*$
and
$ (\h11\h12\h13)^*$ in \kko{genE} yield no element in $H^{13,12}E(\h ij)_4$.
We have cocycles whose leading terms are $x$'s in the next six relations $d(x)=y$ generates the submodule $Z$ of Ker $d$ of rank 6.
% generators in \kko{genE} yield the six linearly independent generators in Ker $d$.
The last six relations with the second and the third ones show that the module $Z$ is in the image of $d$,
%show 
%$$
% \rank \im (d\cln E(\h ij)^{12,12} \to E(\h ij)^{13,12})\ \ge 6,
%$$
 and hence 
\fi
Therefore, no generator in \kko{genE} yields any generator of $H^{13,12}E(\h ij)_4$. \q
 
\if0  % Verification of the above claim
 Put the generators 
 $$
 A_{k,\ell}, B_k, C_k; D, E,F,G;H,I,J
 $$
 Then, the differentials are
 \AL{
 d(X_{k,\ell})&=A_{k,\ell}& d(Y_k)&=B_k+C_k &
 d(B_k)&=a_k&d(D)&=b+c\\
 d(E)&=a_1+a_2& d(F)&=a_2+a_3&d(G)&=a_3+a_0&\\
 d(Z_1)&=H+E+C_1+C_0&d(Z_2)&=I+F+C_2+C_3&d(Z_3)&=J+F+E\\
d(Z_4)&=B_1+B_2+E &d(Z_5)&=H+G+F&d(Z_6)&=%J+G+B_0+B_1&d(Z'_6)&=
G+C_3+C_0\\
 d(H)&=a_0+a_2& d(I)&=0 &d(J)&= a_1+a_3
 }

$\wt E=0$ by $Z_4$. $\wt H=0$ by $Z_1$. $\wt I=\wt F$ by $Z_2$. $\wt J=\wt F$ by $Z_3$.
$\wt F=\wt G$ by $Z_5$. $\wt G=0$ by $Z_6$. 
\fi %end of the verification
\if0
 Then, $F+C_2+C_3$ is a cycle and $d(F)=a_2+a_3$ and $d(C_k)=-d(B_k)=-a_k$ for $k=2,3$.
 Thus $I\sim F+C_2+C_3=\wt F$.
 $\wt E=E\pm C_1\pm C_2$ is a cycle. $\wt H=H\pm C_0\pm C_2$ is a cycle and $H+E+C_1+C_0$ is a cycle.
 So $C_2$ has opposite sign. So $ \wt E+\wt H$ is the image of $d(Z_1)$.
 \fi

 \lem{int}{Let $p$ and $n$ satisfy $n^2\le 2p-1$.
Consider an integer $a=\sum_{i=0}^{n-1}\ep_ie(i)\in \Z/e(n)$ $(\ep_i\in\{0,1\})$. 
Then, $-a=n+\sum_{i=0}^{n-1}(p-1-\ep_i)e(i)$.
It also follows that
$a=\sum_{i=0}^{n-2}a_ip^i$ for integers $a_i$ satisfying $0=a_{n-1}\le a_{n-2}\le \dots\le a_1\le a_0=\sum_{i=1}^{n-1}\ep_i$ and $a_{i-1}-a_i=0$ or $1$, and 
 $-a= \sum_{i=0}^{n-2}(p-a_i)p^i+1$.
 }
 \p
 $-a=\sum_{i=0}^{n-1}p^i-\sum_{i=0}^{n-1}\ep_ie(i)=\sum_{i=0}^{n-1}(p^i-\ep_ie(i))$.
 Since $e(i)=\frac{p^i-1}{p-1}$, we have $p^i-\ep_ie(i)=(p-1)e(i)+1-\ep_ie(i)=(p-1-\ep_i)e(i)+1$. %$p^i-e(i)=(p-1)e(i)+1-e(i)=(p-2)e(i)+1$.
 
We verify the latter half as follows: 
% $a=\sum_{i=0}^{n-2}a_ip^i$ with $0\le a_{n-2}\le a_{n-3}\le \dots\le a_0=\sum_{i=0}^{n-1}\ep_i$ and $a_{i-1}-a_i\in\{0,1\}$.
 %Then, 
 $-a=e(n)-a=%p^{n-1}+e(n-1)-a=\sum_{i=0}^{n-2}(p-1)p^i+1+\sum_{i=0}^{n-2}(1-a_i)p^i=
 \sum_{i=0}^{n-2}(p-a_i)p^i+1$.
 \q

\lem{Lan}{Suppose \kko{cond} and consider an integer $a$ in Lemma \ref{int}.
Let $s_\ep=q+2-a_0-\ep$ %and $s'_\ep=q+1-a_0-\ep$
 for $\ep\in\{0,1\}$.
Then, $E(\h ij)_n^{s_\ep,(a+1)q}=0$.
Furthermore, $E(\h ij)_n^{s_\ep-1,(a+1)q}=0$ if \kko{pn} holds.
\if0
Then, $E(\h ij)_n^{s_\ep,(a+1)q}=0$ and $E(\h ij)_n^{s'_\ep,(a+1)q}=0$, where $s_\ep=q+2-a_0-\ep$ and $s_\ep=q+1-a_0-\ep$ for $\ep\in\{0,1\}$
, except for the case $(n,p,\ep)=(3,5,2)$.
\fi
}
\p
Since the exterior algebra $E(\h ij)_n$ has $n^2$ generators, 
$$
E(\h ij)_n^{s,*}=0 \mbx{if $s>n^2$.} \lnr{dim}
$$
Let $\O s_\ep=q+2-a_0-\ep$ for $\ep\in\{0,1,2\}$ with $\ep\ne 2$ unless \kko{pn} holds. %if $n$ is odd, or if $p\cg 1$ $(4)$.

For $n$ with $n^2+n\le 2p-1$, we see that $n^2+n\le 2p-2$ since $n^2+n$ is even, and so
 the lemma follows from $\O s_\ep=q+2-a_0-\ep\ge 2p-n-1>n^2$ by Lemma \ref{int} and \kko{dim}. %this holds by \cite{sy}.
So we put $n^2+u=2p-2=q$ for $-1\le u\le n-1$, and we have $\O s_\ep=n^2+u+2-a_0-\ep$ %$q+2-a_0-\ep=n^2+u-a_0-\ep$.
%

%Note that
By the duality \kko{isoD}, $E(\h ij)_n^{\O s_\ep,(a+1)q}\cong E(\h ij)_n^{a_0-u-2+\ep,-(a+1)q}$.
Note that $a_0-u-2+\ep\le n< p$.
\if0
We first consider the case where $a_0-u+\ep=p$.
Then, $a_0=n-1$, $u=-1$, $\ep=2$ and $n=p-2$, and we obtain $p=5$, $n=3$, $a_0=2$, and so $s_2=6$, and $-(a+1)=3+(5-a_1)5$
$a=\ep_1+\ep_2(p+1)=\ep_1+\ep_2+\ep_2p$, $a_1=\ep_2$, $a_0=\ep_1+\ep_2=2$. $-(a+1)=23$.
$a=7$, $e(3)=5^2+5+1=31$. 
\fi
%So $a_0-u+\ep<p$ in our case.
For the case $a\ne 0$, let $\di=\max\{i\mid a_i\ne 0\}$. Then, $a_\di=1$ and $p-a_{\di}=p-1$, and
 Lemma \ref{int} implies that the dimension of an element of the form $\prod \h ij$ with degree $-(a+1)q$ is at least $p-1$, which is greater than $n\ge a_0-2-u+\ep$ by \kko{cond}. 
\if0
Since $n^2=q+1$,
$E(\h ij)^{q+2-a_0-\ep,(a+1)q}\cong E(\h ij)^{a_0-1+\ep,-(a+1)q}$.
If $a\ne 0$, then let $\di=\max\{i\mid a_i\ne 0\}$. Then, $a_\di=1$ and $p-a_{\di}=p-1$.
Then, Lemma \ref{int} implies that the dimension of an element of the form $\prod \h ij$ with degree $-(a+1)q$ is at least $p-1$, which is greater than $a_0-1+\ep$. 
\fi
Suppose that $a=0$. Then, $a_0=0$ and $\O s_\ep=q+2-a_0-\ep\ge n^2$. % if $\ep=0,1$.
%If $\ep=2$, then $n$ is even or $p\not \cg 1$ $(4)$, and so $u=q-n^2\ge 0$, and so we also have $\O s\ge q\ge n^2$. 
Therefore, the lemma follows from \kko{dim} if $\O s_\ep\ne n^2$. %Since the exterior algebra has $n^2$ generators, $E(\h ij)^{s_\ep}=0$ if $s_\ep>n^2$.

If $\O s_\ep=n^2$, then
 $E(\h ij)_n^{n^2,q}=0$ by \kko{isoD}. %this is the Morava vanishing theorem.
 %Similarly, if $s_\ep=n^2-1$, then  $E(\h ij)_n^{n^2-1,q}=E(\h ij)_n^{1,-q}=E(\h ij)_n^{1,pe(n-1)q}=\Z/p\{h_{n-1,1}\}$
\q

This lemma and Lemma \ref{HS} imply the following:

\cor{Lanc}
{Suppose \kko{cond} and consider an integer $a$ in Lemma \ref{int}.
Then, for $\ep\in\{0,1\}$, $H^{q+2-a_0-\ep,(a+1)q}S(n)=0$,
and $H^{q+1-a_0-\ep,(a+1)q}S(n)=0$ if \kko{pn} holds. %$n$ is even, or unless $p\cg 1$ $(4)$.
In particular, for $(a,\ep)=(0,1)$, we have $H^{q+1,q}S(n)=0$.
}

\section{Recollection from \cite{sy}}

\if0
Let $J=(a_0,\dots,a_{n-1})$ be an invariant regular ideal of $BP_*$, and put
$$
\La_J=\{\sum_{i=0}^{n-1}\ep_i(|a_i|+1)\mid \ep_i\in\{0,1\}\}.
$$
Here $|a|$ denotes the degree of an element $a\in BP_*$, which is given by
$$
|v_i|=2(p^i-1)=e(i)q. 
$$.
\fi
Let $I_k$ denote the invariant prime ideal of $BP_*$ generated by $p, v_1,\dots, v_{k-1}$.
$$
\La_k=\{\sum_{i=0}^{k-1}\ep_i(e(i)q+1)\mid \ep_i\in\{0,1\}\}.
$$
Here, we notice that the degree of $v_i$ is given by
$$
|v_i|=2(p^i-1)=e(i)q. 
$$.

\if0
We consider
the following lemma to show the existence of $L_nV(n-1)$ for the case where $n^2=2p-1$.
In this case, $n$ is odd, and for $n\le 3$, $V(n-1)$ exists.

For $n=5$, $L_5V(4)$ exists for $p=13$.
$L_7V(6)$ exists for $p=29$.
$L_9V(8)$ for $p=41$?

$L_6V(5)$, $6^2+6=42<2p$ holds for $p\ge 23$.
If $p=19$, $6^2=36<2p-1$ ?
$L_8V(7)$, $8^2+8=72<2p$ holds for $p\ge 37$,
If $p=31$, $64>62$.

The following table shows the smallest prime $p$ for each $n$ satisfying $n^2\le 2p-1$:
\AR{|c||c|c|c|c|c|c|c|c|c|c|}{
\hline
n& 4&5&6&7&8&9&10&11&12&13\\
\hline
n^2&16&25&36&49&64&81&100&121&144&169 \\
\hline
p&11&13&19&29&37&41&53&61&73&87\\
\hline
}
In \cite[Th.~5.7]{sy}, it is shown that $L_nV(n-1)$ exists if $n^2+n<2p$.
This together with the above table shows that $L_4V(3)$, $L_7V(6)$ and $L_8V(7)$ exist.

For $n=5$, we have $L_5V(3)$  which is furthermore a ring spectrum. 
By Lemma \ref{zero}, $v_4\in E_2^0(V(3))_5$ is a permanent cycle, and we obtain $L_5V(4)$ as a cofiber of $v_4\cln L_5V(3)\to L_5V(3)$.

So we consider the case for $n\ge 5$. In this case, $n\le p-8$. 
\fi

%For a positive integer $n$, the relation $n^2\le 2p-1$ implies $n<p$, and so $(p-1)\mid n$ if and only if $n =p-1$.

\lem{Ext}{
Suppose that %a prime number $p$ and a positive integer $n$ satisfy 
$n^2\le  2p-1$  and $n\ne p-1$ (\cf \kko{cond}). %Then, $n$ is odd, and put $n=2k+1$. We then have $p=2k(k+1)+1$ 
Then,
$$
\e_{BP_*BP}^{m+2,m+t}(BP_*/I_n, v_n^{-1}BP_*/I_n)=0\mbx{for $m\ge 1$ and $t\in\La_{n}$.}
$$
Furthermore, if \kko{pn} holds, then %$n$ is even, or unless $p\cg 1$ $(4)$, then
$$
\e_{BP_*BP}^{m+1,m+t}(BP_*/I_n, v_n^{-1}BP_*/I_n)=0\mbx{for $m\ge 1$ and $t\in\La_{n}$.}
$$
\if0
If $n^2<2p-1$, then this holds for any invariant regular ideal $J_n$ of length $n$:
$$
\e_{BP_*BP}^{m+2,m+t}(BP_*/J_n, v_n^{-1}BP_*/J_n)=0\mbx{for $m\ge 1$ and $t\in\La_{J_n}$.}
$$
\fi
}

\p
\if0
$m+t=kq$. then $m=kq-t=kq-aq-b$ with $0\le b\le n$.
Since $m\le 2p-3$, $0\le (k-a)q-b\le 2p-3$. $k-a=1$.
$m=q-b$, $q-n\le m\le q$.
Then the result follows unless $m+b=q$. If $m+b=q$, then $b=0$ and $m=q$. Ext$^{q+2,q}$  This is zero since $q+2>2p-1\ge n^2$. 

$a=0,1,p+1, p^2+p+1,$
\fi
Put $H^{s,t}M=\e_{BP_*BP}^{s,t}(BP_*,M)$ for a $BP_*BP$-comodule $M$, and $A_n=v_n^{-1}BP_*/I_n$.
Then, $\e_{BP_*BP}^{s,t}(BP_*/I_n, v_n^{-1}BP_*/I_n)\cong H^{s,t}A_n$, which is isomorphic to $\LR{K(n)_*\ox H^{*,*}S(n)}^{s,t}$ 
by the Miller-Ravenel change of rings theorem.
%As noted above, $n^2=2p-1$ implies $n\le p-2$.
Since $(p-1)\nmid n$ by \kko{cond}, %It follows that 
$H^{m+e,*}A_n=0$ for $e\in \{1,2\}$ if $m+e> n^2$ by the Morava vanishing theorem.
We assume that $m+e\le n^2$.
By the sparseness, we see that $H^{m+e,m+t}A_n=0$ unless $m+t\cg 0$ mod $q$.
We further assume that $m+t=kq$ for some integer $k$.
%
\if0
$n^2+n-k\le 2p-1=q+1$.
$1\le m\le n^2-2\le 2p-3+k-n\le 2p-3$, $m=q-b$ and $m+t=(a+1)q$.
$$
m+2=q+2-b\ge n^2+n-k+1-b\ge n^2-k+1
$$
So we consider $k=3$. % If $n^2=2p-1$, $n=\sqrt{2p-1}=p-x$. $2p-1=p^2-2px+x^2$ $x=p\pm\sqrt{p^2-(p-1)^2}=p\pm\sqrt{2p-1}$ 
$n^2+n-3\le q+1$.
$m+2=q+2-b\ge n^2+n-3+1-b>n^2$ if $b<n-2$, So we consider $b=n-2$, $n-1$ and $n$.
\fi
%
%
Put $t=aq+b$ with $0\le b<q$. Note that $a$ is an integer of Lemma \ref{int}, and $0\le b\le n\ (\le p-2$ by \kko{cond}$)$ by the definition of $\La_{n}$. Then, $m=(k-a)q-b$.
Since $1\le m\le n^2-2\le 2p-3$, we see $m=q-b$, and so $m+t=(a+1)q$. % and $q-b+2=m+2\le n^2=q+1$.
\if0
Since $1\le m\le 2p-3$, we see $m=q-b$, and so $m+t=(a+1)q$ and $q-b+2=m+2\le n^2=q+1$.
Furthermore, we note that $a\cg b, b-1$ mod $p$, since $|v_i|=e(i)q\in \Z/(e(n)q)$.% and $b<p$.
\fi
Furthermore, let $a_0$ denote an integer such that $a\cg a_0$ mod $p$ in Lemma \ref{int}, and we see that
 $b-a_0=\sum_{i=0}^{n-1}\ep_i-\sum_{i=1}^{n-1}\ep_i=\ep_0$ by the definition of $\La_{n}$.
Thus, $m+e=q+e-b=q+e-a_0-\ep_0$.
Therefore, Corollary \ref{Lanc} implies the lemma. 
\q

\noindent
{\it Proof of Theorem \ref{main3}.}
%$n^2\le 2p-1\le n^2+n-1$ implies $\frac{n^2+1}2\le p\le \frac{n^2+n}2$, $1\le \frac{2p}{n^2+1}\le \frac{n^2+n}{n^2+1}=1+\frac{n-1}2$.
Under Lemma \ref{Ext}, a similar proof of 
 \cite[Th.~5.7]{sy} (using \cite[Th.~4.6]{sy} instead of \cite[Th.~4.9]{sy}) works to show the theorem. \q %  \ref{main3}. \q%the following theorem, which is Theorem \ref{main3}:

\if0
\thm{sy46}{Let $n$ be a positive integer satisfying $n^2\le 2p-1$. %, and $J$ denote an invariant regular ideal of length $n$ in $BP_*$..
Then, there exists a spectrum $L_nV(n-1)$.
%Furthermore, there exists a $v_n$-map $ \Si^{s|v_n|}L_nV(n-1)\to L_nV(n-1)$ for each $s\in \Z$, which is a homotopy equivalence. 
%
%If $n^2<2p-1$, then $L_nMJ_n$ exists, where $J_n$ is an invariant regular ideal of $BP_*$.
}
\fi

\section{The dual map $D_n(i_n)\cln D_n(W_n)\to S$}

In this section, we assume \kko{cond}. % that $n^2\le 2p-1$ and $(p,n)\ne (2,1)$, $(3,2)$.
Then by Theorem \ref{main3}{, we have the spectrum $W_n$.
Let $D_n(-)=F(-,L_nS^0)$ denote the Spanier-Whitehead duality functor on the subcategory $\cL_n$ of $\cS$ consisting of $E(n)$-local spectra.
By \cite[Th.~5.11]{heard}, we see that $W_n$ is strongly dualizable.
Then, we have maps
$$
e\cln D_n(W_n)\sm W_n\to S\qand c\cln S\to W_n\sm D_n(W_n) 
$$
such that the composites 
\AR c{
W_n\xar{c\sm W_n}W_n\sm D_n(W_n)\sm W_n\xar{W_n\sm e} W_n\qand\\
 D_n(W_n)\xar{D_n(W_n)\sm c}D_n(W_n)\sm W_n\sm D_n(W_n)\xar{e\sm D_n(W_n)}D_n(W_n)
}
are the identities on $W_n$ and $D_n(W_n)$, respectively.

Note that the $E_2$-term of the $E(n)$-based Adams spectral sequence for computing $\pi_*(W_n)$ is
$$
E_2^{s,t}(W_n)_n\cong K(n)_*\ox H^{s,t}S(n) \lnr{E2W}
$$
since $E(n)_*(W_n)=E(n)_*/I_n=K(n)_*$.

\lem{ext}{ Let $J=(p^{e_0},v_1^{e_1},\dots, v_{n-1}^{e_{n-1}})$  be an invariant regular ideal of $BP_*$.
Suppose that there exists a spectrum $WJ\in \cL_n$ such that $BP_*(WJ)=v_n^{-1}BP_*/J$.
Then, 
$E(n)_*(D_n(WJ))\cong g_JE(n)_*/J$ for a generator $g_J$ with $|g_J|=n+\sum_{i=0}^{n-1}|v_i^{e_i}|$.
}
\p
Consider  the Adams universal coefficient theorem \cite[Th.~13.6]{ad}
$$
\e_{BP_*}^{s,*}(BP_*(X), E(n)_*)\Longrightarrow E(n)^{s+*}(X)
$$
for the ring spectrum $BP$.
The spectrum $WJ$ is strongly dualizable in $\cL_n$ by \cite[Th.~5.11]{heard}, and then $E(n)^*(WJ)=E(n)_*(D_n(WJ))$.
Put
$$
\Ex^{s,t}(M)=\e_{BP_*}^{s,t}(M,E(n)_*) \mbx{for a $BP_*$-module $M$.}
$$
We notice that $\Ex^{s,t}(v_n^{-1}BP_*/J)\cong \Ex^{s,t}(BP_*/J)$, and
show
%By the Adams universal coefficient theorem \cite[Th.~13.6]{ad}, it suffices to show the following:
$$
\Ex^{s,t}(BP_*/J_k)=\cass{(E(n)_*/J_k)_{t+a_k}&s=k\\0&\text{otherwise}}
$$
by induction on $k\le n$. Here $J_k$ is an ideal of $BP_*$ and $a_k$ is an integer defined by
$$
J_k=(p^{e_0},v_1^{e_1},\dots, v_{k-1}^{e_{k-1}})\qand  a_k=\sum_{i=0}^{k-1}|v_i^{e_i}|.
$$
For $k=0$, it is trivial, since $BP_*$ is a projective $BP_*$-module.
Suppose that it holds for $k$.
Then, we verify the equality for $k+1$
by the long exact sequence
\AL{
\cdots \to \Ex^{s,t}(BP_*/J_k)&\xar{(v_k^{e_k})^*}\Ex^{s,t+|v_k^{e_k}|}(BP_*/J_k) \xar{\de}\Ex^{s+1,t}(BP_*/J_{k+1})\\
&\to \Ex^{s+1,t}(BP_*/J_k)\xar{(v_k^{e_k})^*}\Ex^{s+1, t+|v_k^{e_k}|}(BP_*/J_k)\to\cdots
}
associated to the short exact sequence $0\to BP_*/J_k\xar{v_k^{e_k}}BP_*/J_k\to BP_*/J_{k+1}\to 0$ of $BP_*$-modules.
\if0

$$
\e_{BP_*}^{s,t}(BP_*/I_n,E(n)_*)=\cass{E(n)_*/I_n&(s,t)=(n,-\sum |v_i|),\\0&\text{otherwise}},
$$
since $\e_{BP_*}^{s,t}(v_n^{-1}BP_*/I_n,E(n)_*)\cong \e_{BP_*}^{s,t}(BP_*/I_n,E(n)_*)$.

Put $H^{s}M=\e_{BP_*}^{s}(M,E(n)_*)$ for a $BP_*$-module $M$.
Inductively, we show that
$$
H^sBP_*/I_k=\cass{E(n)_*/I_k&s=k,\\0&\text{otherwise}}.
$$
This holds for $k=0$, since $BP_*$ is a projective $BP_*$-module.
Suppose that it holds for $k$.
Now the lemma follows from the long exact sequence
\AL{
\cdots \to H^{s,t}BP_*/I_k&\xar{v_k^*}H^{s,t+|v_k|}BP_*/I_k \xar{\de}H^{s+1,t}BP_*/I_{k+1}\\
&\to H^{s+1,t}BP_*/I_k\xar{v_k}H^{s+1, t+|v_k|}BP_*/I_k\to\cdots
}
associated to the short exact sequence $0\to BP_*/I_k\xar{v_k}BP_*/I_k\to BP_*/I_{k+1}\to 0$ of $BP_*$-modules.
\fi
\q

\lem{E2ex}{Suppose \kko{cond}. % that $n^2\le 2p-1$ and $n\le p-2$.
Then, $E_2^{s+q, tq}(W_n)_n= 0$ unless 
\begin{itemize}
\item $(s,n^2)=(1,q+1)$ or $(0,q)$, and $t\cg 0$  $(e(n))$.
\item $(s,n^2)=(0,q+1)$, and $t\cg -p^j$  $(e(n))$ for $0\le j<n$.
\end{itemize}
}

\p
By \kko{E2W} and Lemma \ref{HS}, it suffices to show $H^{s+q, tq}E(\h ij)_n=0$ unless the conditions in the lemma.

If $s\ge 2$, then $s+q\ge 2p>n^2$, and so $H^{s+q, tq}E(\h ij)_n=0$.
Suppose that $s=1$ and $n^2\le q$, or that $s=0$ and $n^2\le q-1$. Then $q+s>n^2$, and $H^{s+q, tq}E(\h ij)_n=0$.
For $(s,n^2)=(1,q+1)$ or $(0,q)$, $H^{q+s,*}E(\h ij)_n\cong \Z/p\lrk{g_n}$. Since $|g_n|\cg 0$ mod $2(p^n-1)$, the lemma holds in this case.
For $(s,n^2)=(0, q+1)$, $E(\h ij)_n^{q,tq}\cong \Z/p\lrk{\h ij^*}$.
Since $d((\h 1j\h{i-1}{1+j})^*)=\pm \h ij^*$, we see $\h ij^* =0\in H^{q,tq}E(\h ij)_n$ unless $i=1$.
It follows that $H^{q,-p^jq}E(\h ij)_n\cong \Z/p\lrk{\h 1j^*}$.
\q

\lem{htpy}{
Let $J_k=(p^{e_0},v_1^{e_1},\dots, v_{k-1}^{e_{k-1}})$ for $0\le k\le n$ be an invariant regular ideal of $E(n)_*$, and suppose that the generalized Moore spectrum $MJ_k$ exists.
Then, we have an isomorphism
$$
[MJ_k, W_n]_{0}\cong %E_2^{s,t}(W_n\sm D_n(MJ_n))=
E_2^{0,0}(W_n)_n\o+\Op_{u(\ne 0)\in \La(J_k)} E_2^{q-s(u),\O u+q}(W_n)_n
$$
under \kko{cond}. % if $n^2\le 2p-1$.
\if0
The $E(n)$-based Adams spectral sequence shows an isomorphism of
the homotopy groups $[MJ_n, W_n]_{t-s}\cong E_2^{s,t}(W_n\sm D_n(MJ_n))=E_2^{s,t}(W_n)\o+\Op_{u(\ne 0)\in \La(J_n)} E_2^{s+q-s(u),t+\O u+q}(W_n)$
\fi
Here, 
$$
\La(J_k)=\left\{\sum_{i=0}^{k-1}\ep_i(e_i(2p^i-2)+1)\ \middle|\  \ep_i\in\{0,1\}\right\},
$$ 
and for $u=\sum_{i=0}^{k-1}\ep_i(e_i(2p^i-2)+1)\in\La(J_k)$,
$$
s(u)=\sum_{i=0}^{k-1} \ep_i\qand \O u=u-s(u).
$$.
}

\p Consider the spectral sequence \kko{ANSS} for $X=W_n\sm D(MJ_k)$:
$$
E_2^{s,t}=\e_{E(n)_*E(n)}^{s,t}(E(n)_*, E(n)_*(W_n\sm D(MJ_k)))\Rightarrow [MJ_k,W_n]_{t-s}.
$$
Since
$E(n)_*(W_n\sm D(MJ_k))%=(E(n)_*/I_n)_*(D(MJ_k))~+
=K(n)_*(D(MJ_k))$, we see that
%Thus, inductively 
$$
E(n)_*(W_n\sm D(MJ_k))=\Op_{u\in \La(J_k)} K(n)_{*+u}.
$$
Note that $E_2^{s,t}(W_n)_n=\e_{E(n)_*E(n)}^{s,t}(E(n)_*,K(n)_*)=\LR{K(n)_*\ox H^*S(n)}^{s,t}$. 
Then, 
\AL{
E_2^{s,s}&=\Op_{u\in \La(J_k)} E_2^{s,s+u}(W_n)_n=\cass{E_2^{0,0}(W_n)_n& s=0\\\displaystyle \Op_{\scriptsize \mx{u(\ne 0)\in \La(J_k)\\s+s(u)=q}}E_2^{q-s(u), q+u-s(u)}(W_n)_n&s=q-s(u)\\0&\text{otherwise}}
%&=E_2^{0,0}(W_n)_n\o+ E_2^{q,q}(W_n)_n\o+ \Op_{u(\ne 0)\in \La(J_k)} E_2^{q-s(u),q+u-s(u)}(W_n)_n.
}
by Lemma \ref{E2ex}.
\if0
Indeed, $(s,s+u)$ satisfies $s\le n^2\le q+1$ and $s+u\cg 0$ mod $(q)$.
Then,  
$$
\Op_{s\ge 0}E_2^{sq,sq}=E_2^{0,0}\o+ E_2^{q,q}=E_2^{0,0}(W_n)_n\o+ \Op
$$
Then, $E_2^{0,0}=E_2^{0,0}(W_n)_n$ and $E_2^{q,q}=E_2^{q,q+u}(W)$
$$
E_2^{s,s}=\Op_u E_2^{s,s+u}(W_n)_n=
$$
$s+u=tq$ then $s+\O u+s(u)=tq$ $s+s(u)\cg 0$ $(q)$.

The lemma now follows from Lemmas \ref{ext} and \ref{E2ex} and the fact that $q-s(u)\ge p$. 
\fi
\q
\if0
$$
[MJ_n, W_n]_*\cong 
$$
\p
Put $a_i=v_i^{e_i}$ so that $J_k=(a_0,a_1,\dots,a_{k-1})$. Then $|a_i|=e_i(2p^i-2)$.
 We show the isomorphism of the lemma and  isomorphisms
$$
[MJ_k, W_n]_{w}\cong %E_2^{s,t}(W_n\sm D_n(MJ_n))=
E_2^{q-s(w),\O w+q}(W_n)\o+\Op_{u(\ne 0)\in \La(J_k)} E_2^{q-s(w)-s(u),\O w+q+\O u}(W_n)
$$
for $w(\ne 0)\in \La(J_{n})$ 
by  induction on $k$.
The spectra $MJ_k$'s fits in cofiber sequences 
$$
\Si^{|\ph_k|}MJ_k\xar{\ph_k}MJ_k\to MJ_{k+1}
\lnr{cofMJ}
$$
 for $v_k$-map $\ph_k\cln \Si^{|\ph_k|}MJ_k\to MJ_k$ with $E(n)_*(\ph_k)=a_k$. % v_k^{e_k}$.
%Note that $E_2^{s,t}(W_n)\cong K(n)_*\ox H^{s,t}S(n)$, which is zero if $s\ge q+1\ge n^2$.
For $k=0$, $MJ_0=S^0$, and $[S^0,W_n]_{0}\cong E_2^{0,0}(W_n)$ by Lemma \ref{E2ex}
%the lemma holds, since the $E(n)$-based Adams spectral sequence collapses from the $E_2$-term.
and $[S^0,W_n]_{1}\cong E_2^{q-1,q}(W_n)$.

Suppose that the isomorphisms for $k$ hold and consider an exact sequence
$$
[MJ_k, W_n]_{u+1}\xar{\ph_k^*} [MJ_k, W_n]_{u+1+|\ph_k|}\xar{} [MJ_{k+1}, W_n]_{u}\xar{}[MJ_k, W_n]_{u}\xar{\ph_k} [MJ_k, W_n]_{u+|\ph_k|}
$$
associated to the cofiber sequence \kko{cofMJ}.
It follows that $\ph_k^*=0$ and we obtain the desired isomorphism on $MJ_{k+1}$. \q
\fi

\if0
We see that 
$$
E(n)_*(W_n\sm D_n(MJ_n))\cong K(n)_*\ox E(a_k\mid 0\le k< n)
$$
as an $E(n)_*E(n)$-comodule, where $a_i$ is a primitive element with  $|a_k|=e_k|v_k|+1$.
Furthermore, $W_n$ is an $M(p)$-module spectrum.
Therefore, the lemma follows. \q
\fi

Since $E_2^{s,*}(W_n)_n=0$ if $s\ge n^2$, 
this lemma implies
$$
[MJ_n,W_n]_0=E_2^{0,0}(W_n)_n\o+\Op_{u(\ne 0)\in \La(J_n)}E_2^{q-s(u), \O u+q}(W_n)_n.
$$
Then, we have a map $\ph\cln MJ_n\to W_n$ detected by $1\in E_2^{0,0}(W_n)_n$, which satisfies 
$$
\ph i_J=i_n \lnr{comm}
$$ for the maps $i_J\cln S^0\to MJ_n$ and $i_n\cln S^0\to W_n$ detected by $1\in E_2^{0,0}(MJ_n)$ and $1\in  E_2^{0,0}(W_n)_n$, respectively. 

Lemma \ref{ext} for $J=I_n$ shows that $E(n)_*(\wt W_n)\cong E(n)_*/I_n$ for $\wt W_n=\Si^{d_{I_n}}D_n(W_n)$.
Here, $d_{I_n}=2(e(n)-n)+n$.
We have a finite spectrum $MJ_n$ such that   $D_n(MJ_n)=\Si ^{-d_{J_n}}L_nMJ_n$ for $d_{J_n}=n+\sum_{i=1}^{n-1}|v_i^{e_i}|$.
We write $\ph^*\cln \Si ^{d_{J_n}-d_{I_n}}\wt W_n\to L_nMJ_n$ for the map $D_n(\ph)\cln D_n(W_n)\to D_n(MJ_n)$.

\lem{ph}{Put $V_{J_n}=p^{e_0-1}v_1^{e_1-1}\dots v_{n-1}^{e_{n-1}-1}\in E_2^{0,*}(MJ_n)_n$. Then, 
the map $\ph^*\cln \Si ^{d_{J_n}-d_{I_n}}\wt W_n\to L_nMJ_n$ induces $(\ph^*)_*=uV_{J_n}\cln E_2^{0,*}(\wt W_n)_n\to E_2^{0,*}(MJ_n)_n$ for some $u\in\Z/p\setminus \{0\}$. 
}

\p
By Lemma \ref{ext}, the $E(n)_*(-)$-homology induces a non-trivial homomorpshim
$
(\ph^*)_*\cln E(n)_*/I_n\to E(n)_{*+|V_{J_n}|}/J_n$ ($|V_{J_n}|=d_{J_n}-d_{I_n}$)
of $E(n)_*$-modules.
Thus, the lemma follows, 
since $p^{e_0-1}E(n)_{|V_{J_n}|}/J_n\cong \Z/p\{V_{J_n}\}$. \q

\section{The Greek letter elements}
For the invarinat regular ideal $I_k$,
 we have short exact sequences
$$
0\to BP_*/I_k\xar {v_k} BP_*/I_k\to BP_*/I_{k+1}\to 0
$$
of $BP_*BP$-comodules
for $k\ge 0$ ($v_0=p$ and $I_0=(0)$).
Let $\der_k\cln H^sBP_*/I_{k+1}\to H^{s+1}BP_*/I_{k}$ be the connecting homomorphism associated to the above short exact sequence. 
Here, 
$$
H^sM=\e_{BP_*BP}^s(BP_*, M)
$$
for a $BP_*BP$-comodule $M$, and so $H^sBP_*(X)=E_2^s(X)$, the $E_2$-term of the \ANSS\ for $\pi_*(X)$.
Then, we define the $n$-th Greek letter elements $\al^{(n)}_s\in H^nBP_*=E_2^{n,*}(S^0)$ by
$$
\al^{(n)}_s=\der_0\der_1\cdots\der_{n-1}(v_n^s)\mbx{for $v_n^s\in H^0BP_*/I_n$.}
$$
Similarly, take an invariant regular ideal $J_k=(p^{e_0},v_1^{e_1},\dots, v_{k-1}^{e_{k-1}})$ instead of $I_k$.
For the connecting homomorphisms $\der^J_k\cln H^sBP_*/J_{k+1}\to H^{s+1}BP_*/J_k$, the $n$-th Greek letter element $\al_s^{(n)}$
is also given by
$$
\al^{(n)}_s=\der^J_0\der^J_1\cdots\der^J_{n-1}(V_{J_n}v_n^s)\mbx{for $V_{J_n}v_n^s\in H^0BP_*/J_n$.}
$$
Here, $V_{J_n}=p^{e_0-1}v_1^{e_1-1}\dots v_{n-1}^{e_{n-1}-1}$.

Replacing $BP_*$ with $E(n)_*$,
%Similarly, 
the $n$-th Greek letter elements $\al^{(n)}_s$ are also defined in $E_2^n(S)_n$ for any non-zero integer $s$.
Consider the canonical map $\ell_n\cln BP\to E(n)$, which induces the map form the Adams-Novikov spectral sequence to the $E(n)$-based spectral sequence.
Then, $(\ell_n)_*\cln E_2^{s,t}(S^0)\to E_2^{s,t}(S^0)_n$ assigns an $E(n)$-local Greek letter element with positive suffix $s$ to the Greek letter element with the same suffix $s$.
\if0
, that is, 
 we obtain the $n$-th Greek letter element  $\al^{(n)}_s\in H^nE(n)_*=E_2^{n,*}(L_nS^0)$ as the image of $\al^{(n)}_s\in H^nBP_*$ under  the homomorphism $\ell_n\cln H^nBP_*\to H^nE(n)_*$.
\fi
For $n\le 4$, we have the Smith-Toda spectra $ V(n-1)$ satisfying $E(n)_*(V(n-1))=E(n)_*/I_{n}$ if $p\ge 2n-1$.

We note that the Smith-Toda spectrum $V(3)$ exists at the prime 7, and $W_n$ exists if $n^2\le 2p-1$ by Theorem \ref{main3}. %\ref{sy46}.

\lem{perm}{Let $s\in\Z\setminus \{0\}$. Then, the generators
$v_n^s\in E_2^0(W_n)_n$ for $n^2\le 2p-1$ and $v_4^s\in E_2^0(V(3))_4$ for $p\ge 7$ are permanent cycles.
}

\p
\if0
The lemma holds except for the cases $v_n^s$ for $n^2=p=5$ and $v_4$ for $p=7$.

At $p=5$, $H^{9,8}S(3)=0$
\fi
Since $E_2^*(W_n)_n\cong K(n)_*\ox H^*S(n)$, we see that $d_{q+1}(v_n^s)=0$ if $n^2\le q+1$ by Corollary \ref{Lanc}, or if $(p,n)=(7,4)$ by Lemma \ref{zero}.
Furthermore, $H^sS(n)=0$ for $s>n^2$ by the Morava vanishing theorem, and so $d_r(v_n^s)=0$ for $r>n^2$. 
\q

\noindent
{\it Proof of Theorem \ref{main2}}.
Let $A_s^{(n)}\in \pi_*(W_n)$ be the element detected by the permanent cycle $v_n^s$ in Lemma \ref{perm}. 
Consider an element $A=D_n(i_n)A_s^{(n)}\in \pi_*(L_nS)$.
By \kko{comm}, we have $D_n(i_n)=D_n(i_J)D_n(\ph)$.
Then
$$
A=D_n(i_n)A_s^{(n)}=D_n(i_J)D_n(\ph)A_s^{(n)}.
$$
Hopkins and Smith \cite{hs} show that there exists an invariant regular ideal $J_n$ of length $n$ such that the generalized Moore spectrum $MJ_n$ exists.
We see that $D_n(\ph)A_s^{(n)}\in \pi_*(L_nMJ_n)$ is detected by $uV_{J_n}v_n^s\in E_2^0(MJ_n)$ for some $u\ne0\in \Z/p$ by Lemma \ref{ph}.
Since $D_n(i_J)$ induces the composite $\der^J_0\dots \der^J_{n-1}$ on the $E_2$-terms, the Geometric Boundary Theorem implies that $A$ is detected by $u\al_s^{(n)}$.

\q 
\if0
Let $i_n\cln S^0\to W_n$ denote the map given by $1\in \pi_0(W_n)$. Then, we have its Spanier-Whitehead dual $j=D_n(i_n)\cln D_n(W_n)\to D_n(S^0)=L_nS^0$, where $D_n(X)=F(X,L_nS^0)$.
Lemma \ref{perm} shows that $j_*(v_n^s)\in \pi_*(L_nS^0)$.
Then, by Lemma \ref{conn} , we see that $\al^{(n)}_s$ detects $j_*(v_n^s)$, that is, $\al^{(n)}_s$ is a permanent cycle.
\q
\fi

\noindent
{\it Proof of Theorem \ref{main}}.
This follows from Lemma \ref{perm} and the Geometric Boundary Theorem. \q

\bibliographystyle{amsplain}

\end{document}